\documentclass[11pt]{article}
\usepackage{amsmath, amssymb}
\newtheorem{theorem}{Theorem}[section]
\newtheorem{lemma}[theorem]{Lemma}
\newtheorem{corollary}[theorem]{Corollary}

\begin{document}
\title{Effective classes and Lagrangian tori in symplectic four-manifolds}
\author{Jean-Yves Welschinger\\
}
\maketitle

{\bf Abstract: }

An effective class in a closed symplectic four-manifold $(X, \omega)$ is a 
two-dimensional homology class
which is realized by a $J$-holomorphic cycle for every tamed almost complex structure $J$.
We prove that effective classes are orthogonal to Lagrangian tori in $H_2 (X ; \Bbb{Z})$.

\subsection{Results}

Let $(X, \omega)$ be a closed symplectic four-manifold. A two-dimensional homology class
$d \in H_2 (X ; \Bbb{Z})$ is called an effective class if it is realized by a $J$-holomorphic
two-cycle for every almost-complex structure $J$ tamed by the symplectic form $\omega$.

\begin{theorem}
\label{theo}
Let $(X, \omega)$ be a closed symplectic four-manifold. Let $A_\omega$ be the subspace of $H_2 (X ; \Bbb{Z})$ generated by Lagrangian tori and $B_\omega$ the subspace generated by effective classes.
Then, $A_\omega$ and $B_\omega$ are orthogonal to each other.
\end{theorem}

\begin{corollary}
\label{cor}
Let $L$ be a Lagrangian torus and $S$ be an embedded symplectic $(-1)$-sphere in a closed symplectic four-manifold $(X, \omega)$. Then, L and S have vanishing intersection index.  $\square$
\end{corollary}

We indeed know from Lemma $3.1$ of \cite{McD} that such embedded symplectic $(-1)$-spheres define effective classes. 
Do there exist a Lagrangian torus and symplectic $(-1)$-sphere such that, though they have vanishing intersection index, they have to intersect?
Otherwise, it means that the space $A_\omega$ comes from an underlying minimal symplectic four-manifold.

Let $L$ be a torus equipped with a flat metric, $S^* L$ be its unit cotangent bundle and $\pi : S^* L
\to L$ the canonical projection. The manifold $S^* L$ is equipped with a canonical contact form 
$\lambda$,
namely the restriction of the Liouville one-form of its cotangent bundle. We denote by $R_\lambda$ the
subgroup of $H_1 (S^* L ; \Bbb{Z})$ generated by its closed Reeb orbits.

\begin{lemma}
\label{lemma}
The restriction of $\pi_* : H_1 (S^* L ; \Bbb{Z}) \to H_1 (L ; \Bbb{Z})$ to $R_\lambda$ is an isomorphism.
\end{lemma}

{\bf Proof:}

The Reeb flow on $S^* L$ coincides with the geodesic flow. Closed Reeb orbits are thus the lifts
of closed geodesics on L. Now $S^* L$ is diffeomorphic to a product of $L$ with the sphere $S$ 
of directions in $L$, and $\pi$ is the projection onto the first factor. Since geodesics of $L$ have a
constant direction, the projection onto the second factor maps every Reeb orbit to a point of $S$. From
K\"unneth formula, we get the isomorphism $H_1 (S^* L ; \Bbb{Z}) \cong H_1 (L ; \Bbb{Z}) \times
H_1 (S ; \Bbb{Z})$ and, from what we have just noticed, that this isomorphism maps $R_\lambda$ into
$H_1 (L ; \Bbb{Z}) \times \{ 0 \}$. Since generators of $H_1 (L ; \Bbb{Z})$ can be realized by closed geodesics,
the latter map is onto. Since $\pi_*$ is the projection onto the first factor, it is an isomorphism once
restricted to $R_\lambda$. $\square$ \\

{\bf Proof of Theorem \ref{theo}:}

Let $L$ be a Lagrangian torus and $d$ be an effective class. Following the principle of
symplectic field theory \cite{EGH}, we stretch the neck of the symplectic manifold in the neighbourhood
of $L$ until the manifold splits in two parts, one part being the cotangent bundle of the torus and the
other one being $X \setminus L$. We produce this splitting in such a way that both parts
have the contact manifold $(S^* L , \lambda)$ at infinity. Let $J_\infty$ be a $CR$-structure on this
contact three-manifold which we extend to an almost-complex structure $J$ with cylindrical end
on both parts $T^* L$ and $X \setminus L$. The latter is the limit of a sequence $J_n$ of almost-complex structures of $(X, \omega)$. Since $d$ is effective, we may associate a sequence $C_n$
of $J_n$-holomorphic two-cycles homologous to $d$. From the compactness Theorem in SFT
\cite{BEHWZ}, we extract a subsequence converging to a broken $J$-holomorphic curve $C$, which
we assume for convenience to have only two levels -the general case follows easily from this one-.
Denote by $C^L$ the part of $C$ in $T^* L$ and by $C^X$ the part in $X \setminus L$. Both curves
$C^L$ and $C^X$ have cylindrical ends asymptotic to the same set of closed Reeb orbits with
same multiplicities. Let $C^L_1, \dots , C^L_k$ denote the irreducible components of $C^L$,
and $R_1, \dots , R_k$ be the corresponding set of closed Reeb orbits. These sets  $R_1, \dots , R_k$ 
define integral one-cycles in $S^* L$ and we denote by $[R_1], \dots , [R_k]$ their homology classes. 
These one-cycles are boundaries of the two-chains $C^L_1, \dots , C^L_k$
in $T^* L$ , so that with the notation of Lemma \ref{lemma}, $\pi_* ([R_i])$ vanishes for every $i \in \{ 1 , \dots
, k \}$. Since $[R_1], \dots , [R_k]$ belong to the subgroup $R_\lambda$, we deduce from Lemma \ref{lemma} that 
$[R_1], \dots , [R_k]$ 
actually vanish. Let $S_1, \dots , S_k$ be integral two-chains of  $S^* L$ having $R_1, \dots , R_k$
as boundaries, and $S$ be the sum of these $k$ chains. Then, $C^L - S$ is an integral two-cycle
contained in $T^* L$, $C^X + S$ is an integral two-cycle contained in $X \setminus L$, and the sum of these cycles is
homologous to $d$. Now, $L$ and $C^X + S$ are disjoint from each other and the second homology
group of $T^* L$ is generated by $[L]$ itself. Since the latter has vanishing self-intersection, we
deduce that the intersection index of $L$ and $C^L - S$ vanishes. As a consequence, the
intersection index of $d$ and $[L]$ vanishes. Since this holds for any Lagrangian torus or
effective class, Theorem \ref{theo} is proved. $\square$ 

\subsection{Remarks}
\label{rem}

\begin{enumerate}
\item We have actually proved more, namely for a class $d \in H_2 (X ; \Bbb{Z})$ to be orthogonal to a Lagrangian torus $[L]$, it suffices that $d$ be realized by a sequence of $J_n$-holomorphic two-cycles, for a sequence $J_n$ having a flat neck stretching to infinity.

\item  From the results of Taubes \cite{T}, any Seiberg-Witten basic class is effective and thus,
from Theorem \ref{theo}, we deduce that SW-basic classes are orthogonal to Lagrangian tori.
This fact was already known, it indeed follows from the adjunction inequality \cite{KM}, \cite{MST}. 
Our space $B_\omega$ might  however be bigger than the one generated by SW-basic classes?
Also, our proof remains in the symplectic category and offers possibilities to have counterparts in 
higher dimensions.

\item  If $(X, \omega)$ is K\"ahler, then the Poincar\'e dual of $B_\omega$ is contained in the intersection of $H^{1,1} (X ; \Bbb{Z}) $ on every
complex structure of $X$ tamed by $\omega$. How smaller can it be? 

\item From Hodge' signature Theorem follows that the intersection $A_\omega \cap B_\omega$ vanishes for every K\"ahler surface. 
This intersection indeed lies in the isotropic cone of the Lorentzian $H^{1,1} (X ; \Bbb{R})$ and
is orthogonal to the symplectic form which lies in the positive cone (compare Example $1$ of \S \ref{example}).
I don't see at the moment whether or not this holds for every closed symplectic four-manifold (and am grateful to St\'ephane Lamy for raising
the question to me). More generally, one may wonder whether the intersection form restricted to $B_\omega$ has to be non-degenerated. 
This is the case at least for rational surfaces from Example $2$ of \S \ref{example} and for K\"ahler surfaces with $b_2^+ \geq 2$ and $K_X^2 > 0$ 
since from Taubes' results \cite{T}, $B_\omega$ contains the canonical class $K_X$.

\item  We made a crucial use of a property of the contact manifold $(S^* L , \lambda)$, namely that the
subgroup $R_\lambda$ generated by its closed Reeb orbits has a rather big index in $H_1 (S^* L ; \Bbb{Z})$.
We may then more generally wonder, given a contact manifold, how small can this subgroup $R_\lambda$ 
of  "effective" homology classes be?
\end{enumerate}

\subsection{Examples}
\label{example}

\begin{enumerate}
\item If  $(X, \omega)$ has a Lorentzian intersection form, then $A_\omega$ vanishes. Indeed, each Lagrangian
torus should be in the isotropic cone of the intersection form and should be orthogonal to the class
of the symplectic form which lies in the positive cone.

\item If  $(X, \omega)$ is a blow up of the projective plane, then $B_\omega = H_2 (X ; \Bbb{Z})$. Indeed,
exceptional divisors are effective classes, and the strict transform of a line has non-trivial GW-invariants.

\item If $X$ is a product of two curves $(C_1 , \omega_1)$ and $(C_2, \omega_2)$ with symplectic form
$\omega_1 \ominus \omega_2$, then $A_\omega$ contains the index two subgroup 
$H_1 (C_1 ; \Bbb{Z}) \otimes H_1 (C_2 ; \Bbb{Z})$ given by K\"unneth formula. When in addition,
$(C_1 , \omega_1)$ and $(C_2, \omega_2)$ are symplectomorphic tori,  $A_\omega$ contains the graph of 
the symplectomorphism. Note that $A_\omega$ cannot have codimension less than one since it lies in the orthogonal of the symplectic form.

\item If $X$ is a product of two genus one curves, then $B_\omega$ vanishes, since for a generic complex structure,
$H^{1,1} (X ; \Bbb{Z}) $ vanishes. If instead one of the curve is not elliptic, then we know from Taubes' results  \cite{T} that 
$B_\omega$ contains the canonical class of $X$.
\end{enumerate}

{\bf Acknowledgements:} I am grateful to the French Agence Nationale de la Recherche for its support.

\'Ecole normale sup\'erieure de Lyon\\
Unit\'e de math\'ematiques pures et appliqu\'ees\\
UMR CNRS 5669\\
46, all\'ee d'Italie\\
69364, Lyon cedex 07\\
(FRANCE)\\ 
\end{document}